\documentclass{amsart}

\usepackage{amssymb}
\usepackage[all]{xy}
\usepackage{hyperref}
\usepackage{comment}

\usepackage{enumitem}   %from SK, fixes \enumerate

\setlist[enumerate]{itemsep=.2em,topsep=.2em,leftmargin=1.25em,itemindent=2.0em}

\newtheorem{thm}{Theorem}%[section]

\newtheorem{lem}[thm]{Lemma}
\newtheorem{cor}[thm]{Corollary}

\newtheorem{prop}[thm]{Proposition}

\newtheorem{complem}[thm]{Complement}%!!!!!!!!!!!!!!!!!!!!!!
\theoremstyle{definition}
\newtheorem{defn}[thm]{Definition}

\newtheorem{say}[thm]{}
\newtheorem{exmp}[thm]{Example}

\newtheorem{ques}[thm]{Question}    %!!!!!!!!!!!!!!!!!!!!

\newtheorem*{ack}{Acknowledgments}      % \renewcommand{\theack}{} 

\newtheorem{defn-thm}[thm]{Definition--Theorem}  %!!!!!!!!!!!!!!!!!!!!!!!!
\newtheorem{defn-lem}[thm]{Definition--Lemma}  %!!!!!!!!!!!!!!!!!!!!!!!!
\theoremstyle{remark}

\renewcommand{\c}[0]{{\mathbb C}}  

\renewcommand{\o}[0]{{\mathcal O}} 
\newcommand{\z}[0]{{\mathbb Z}}
\newcommand{\n}[0]{{\mathbb N}}
\renewcommand{\r}[0]{{\mathbb R}}

\newcommand{\dd}[0]{{\mathbb D}}

\newcommand{\p}[0]{{\mathbb P}}

\newcommand{\q}[0]{{\mathbb Q}}

\newcommand{\qtq}[1]{\quad\mbox{#1}\quad}

\newcommand{\mult}[0]{\operatorname{mult}}

\newcommand{\res}[0]{\operatorname{\mathcal R}}

\newcommand{\onto}[0]{\twoheadrightarrow}

\newcommand{\tsum}[0]{\textstyle{\sum}}

\def\into{\DOTSB\lhook\joinrel\to}

\def\loccoh#1.#2.#3.#4.{H^{#1}_{#2}(#3,#4)}

\DeclareMathAlphabet{\mathchanc}{OT1}{pzc}%
                                {m}{it}

\usepackage[all]{xy}\xyoption{dvips}

\begin{document}
\bibliographystyle{amsalpha}
%\hfill\today

  \title{Seshadri's criterion and openness of projectivity}
  \author{J\'anos Koll\'ar}

\begin{abstract} 
We prove that   projectivity is an open condition for deformations of algebraic spaces with rational singularities.
\end{abstract}

 \maketitle

Many projective varieties have deformations that are not projective, not even algebraic in any sense; K3 and elliptic surfaces furnish the best known examples.
In these cases the non-algebraic deformations are also very far from being algebraic.   
A typical non-projective K3 surface does not contain any compact curves and
a non-projective K3 surface  is not  even bimeromorphic to an algebraic surface. The latter property is not an accident.

Let $g: X\to \dd$ be a smooth, proper morphism of a complex manifold to a disc.
If the central fiber $X_0$ is  projective, then it is K\"ahler, and the $X_s$ are also  K\"ahler for $|s|\ll 1$ 
by \cite{MR0112154}. A K\"ahler variety that is bimeromorphic to an algebraic variety is projective by \cite{Moi-66}. Thus if  $X_s$ is bimeromorphic to an algebraic variety, then it  is  projective for $|s|\ll 1$. We can summarize this somewhat imprecisely as:

\medskip

$\bullet$ 
Projectivity is an open condition in the category of Moishezon manifolds.
\medskip

Note, however, that even if all fibers of $g: X\to \dd$ are projective,
$g$ need not be projective over any smaller disc $\dd_\epsilon$; see
\cite{Atiyah58} or \cite[Exmp.4]{k-defgt}.
Examples of Hironaka---reproduced in \cite[App.B]{hartsh}---show that
projectivity is not a closed condition.

The aim of this paper is to prove that   projectivity is an open condition for deformations of algebraic spaces with rational singularities.

A result of this type, with rather strong restrictions on the singularities, is in  \cite[12.2.10]{km-flips}; see Paragraph~\ref{corr.12.2.10} for some comments.
I have been trying to remove the restrictions, and a key problem that emerged was that the argument  used Kleiman's criterion for ampleness, which is not known to apply to algebraic spaces.

By contrast, Seshadri's criterion does hold for   algebraic spaces, and
leads to a quite general openness result.

First I state the 1-parameter version. This is clearer and its proof uses all the essential ideas.
The general version, Theorem~\ref{12.2.10.thm.S}, is substantially stronger.

\begin{thm} \label{12.2.10.thm.dd}
Let $g:X\to \dd$ be a  proper, flat morphism of complex analytic spaces.
Assume that
\begin{enumerate}
  \item $X_0$ is projective, 
\item the fibers $X_s$ have rational singularities for $s\neq 0$, and
\item $g$ is bimeromorphic to a projective morphism $g^{\rm p}:X^{\rm p}\to \dd$.
\end{enumerate}
Then $g$ is projective over a smaller punctured disc  $\dd^\circ_\epsilon\subset \dd$.
\end{thm}

{\it Remarks} (\ref{12.2.10.thm.dd}.4)  Note that we make no assumptions on the singularities of $X_0$.

(\ref{12.2.10.thm.dd}.5)  Even if $g$ is smooth, usually $g$ is not projective over the disc $\dd_\epsilon$. So, although  there is a line bundle $L$ on $X$ such that
$L$ is ample over $\dd^\circ_\epsilon$, it can not be chosen to be also ample on $X_0$. This is quite typical for simultaneous resolutions of deformations of surfaces with Du~Val singularities;  see \cite{Atiyah58} or \cite[Exmp.4]{k-defgt}.

(\ref{12.2.10.thm.dd}.6)  Assumption (\ref{12.2.10.thm.dd}.2) can be relaxed to $X_s$ having 1-rational singularities  (see Definition~\ref{rtl.1.rtl.defn}), but  a similar  openness property does not hold  for  algebraic spaces with slightly worse singularities. 
There are  flat, proper families of normal, compact surfaces 
where all fibers are bimeromorphic to $\p^2$, yet 
the projective fibers correspond to  a countable, dense set on the base; see Example~\ref{bad.fams.exmp}.
In this example each surface has a single singular point, which is simple elliptic (that is, biholomorphic to a cone over an elliptic curve).
These are among the mildest surface singularities, but they are not rational.

(\ref{12.2.10.thm.dd}.7) Assumption (\ref{12.2.10.thm.dd}.3) is automatic if $g:X\to \dd$ is obtained by base change from a morphism of  algebraic spaces,  and it implies that every fiber of $g$ is an algebraic space.
A converse is conjectured in \cite[Conj.2]{k-moi}, and  proved in 
\cite[Thm.1.4]{rao-tsa-1} and \cite[Thm.1.2]{rao-tsa}  for smooth morphisms.

(\ref{12.2.10.thm.dd}.8) Instead of assumption (\ref{12.2.10.thm.dd}.3) we could assume that all fibers are  algebraic spaces and the irreducible components of the Chow-Barlet space parametrizing  1-cycles on $X/\dd$ are proper over $\dd$; see \cite{bar-mag} or  Proposition~\ref{bir.chow.prop}

%% (The latter is somewhat nebulous since there is no construction of this  Chow-Barlet space in the literature. See  Proposition~\ref{bir.chow.prop} for a way to go around this problem.)

\medskip

{\it Idea of proof.} After shrinking $\dd$ we may assume that $X$ retracts to $X_0$. Since $X_0$ is projective,  it has an ample line bundle $L$.
Let  $\Theta\in H^2(X, \q)$ be the pull-back of
$c_1(L)$ to $X$.  Note that $\Theta$ is a topological cohomology class that is usually not the Chern class of a holomorpic line bundle.

Fix a very general $s\in \dd$ and let $p_s\in C_s\subset X_s$ be a pointed curve.
We show in Proposition~\ref{bir.chow.prop} that it has a specialization to $p_0\in C_0\subset X_0$ such that
$\mult_{p_0}C_0\geq \mult_{p_s}C_s$.  Thus
$$
  \Theta\cap [C_s]=\Theta\cap [C_0]
  \geq \epsilon\cdot \mult_{p_0}C_0\geq \epsilon\cdot\mult_{p_s}C_s,
  $$
  where the first $\geq $ is the easy direction of Seshadri's criterion \ref{sesh.crit.thm}  on $X_0$.  Thus $\Theta_s:=\Theta|_{X_s}\in H^2(X_s, \q)$ satisfies the assumption of 
  Seshadri's criterion on $X_s$.

  Thus we need to show that Seshadri's criterion works for (possibly non-algebraic) cohomology classes.  This is done in Proposition~\ref{sesh.top.ver.prop}, using some foundational work comparing numerical and homological equivalence in Proposition~\ref{N1.to.H2.prop}. This is where rationality of the singularities of $X_s$ is used.

  Going from very general fibers to all fibers over a smaller disc  uses Proposition~\ref{RT.vb.Zd.M}.
  
  \medskip

The general version is the following, whose proof is  completed  in Paragraph~\ref{12.2.10.thm.S.pf}.

\begin{thm} \label{12.2.10.thm.S}
  Let $g:X\to S$ be a  proper morphism of complex analytic spaces and
  $S^*\subset S$ a dense, Zariski open subset such that $g$ is flat over $S^*$. 
Assume that
\begin{enumerate}
  \item $X_0$ is projective for some $0\in S$,
\item the fibers $X_s$ have rational singularities for $s\in S^*$, and
\item $g$ is bimeromorphic to a projective morphism $g^{\rm p}:X^{\rm p}\to S$.
\end{enumerate}
Then there is a Zariski  open neighborhood
$0\in U\subset S$ and a locally closed, Zariski stratification
$U\cap S^*=\cup_i S_i$ such that each
$$
g|_{X_i}:X_i:=g^{-1}(S_i)\to S_i \qtq{is projective.}
$$
\end{thm}

{\it Remark \ref{12.2.10.thm.S}.4.}   Note that $g$ is not  assumed flat and
the comments in (\ref{12.2.10.thm.dd}.6) also apply to $X\to S$. 
\medskip

The Lefschetz principle then gives the analogous result for algebraic spaces.

\begin{cor} \label{12.2.10.thm.S.cor} Let $k$ be a field of characteristic 0, and  $g:X\to S$  a  proper morphism of algebraic spaces that are of finite type over $k$. Let 
  $S^*\subset S$ be a dense, Zariski open subset  such that $g$ is flat over $S^*$.
Assume that
\begin{enumerate}
  \item $X_0$ is projective for some $0\in S$, and
\item the fibers $X_s$ have 1-rational singularities for $s\in S^*$.
\end{enumerate}
Then there is a Zariski  open neighborhood
$0\in U\subset S$ and a locally closed, Zariski stratification
$U\cap S^*=\cup_i S_i$ such that each
$$
g|_{X_i}:X_i:=g^{-1}(S_i)\to S_i \qtq{is projective.}\hfill\qed
$$
\end{cor}

\begin{ques} Is Corollary~\ref{12.2.10.thm.S.cor} also true in positive or mixed characteristic? Our proof uses topological cohomology groups in an essential way; it is not clear what would replace them.
  \end{ques}

Theorem~\ref{12.2.10.thm.dd} was developed to complete the proof of the second part of the following; for details see the original paper.

\begin{thm} \cite[Thm.1]{k-defgt} \label{mmp.extends.conj.thm.c.cor}
Let $g:X\to S$ be a flat, proper morphism of complex analytic spaces. 
Fix a point $0\in S$ and assume that the fiber
$X_0$ is  projective, of general type, and with  canonical  singularities.
Then   there is an open  neighborhood  $0\in U\subset S$ such that
\begin{enumerate}
\item the plurigenera  $h^0(X_s, \omega_{X_s}^{[r]})$ are independent of $s\in U$ for every $r$, and
\item  the fibers $X_s$ are  projective  for every $s\in U$.
 \qed
\end{enumerate}
\end{thm}

\begin{exmp} \label{bad.fams.exmp}
Let $E\subset \p^2$ be a smooth cubic. Fix $m\geq 10$ and let  $X\to S$ be the universal family of surfaces obtained by blowing up $m$ distinct points $p_i\in E$, and then contracting the birational transform of $E$.  (So $S$ is an open subset in $E^m$.) Such a surface is projective iff there are positive $n_i$ such that
$\sum_i n_i[p_i]\sim  nL|_E$  where $L$ is the line class on $\p^2$ and $n=\frac13 \sum_i n_i$. 

Thus  the projective fibers 
correspond to a countable union of hypersurfaces
$H(n_1,\dots, n_m)\subset S$.  All fibers have simple elliptic singularities and trivial canonical class.
(For $m=12$ the singularities are biholomorphic to cones over smooth plane cubics.) 
\end{exmp}

\begin{say}[Correction to {\cite[Chap.12]{km-flips}}]
\label{corr.12.2.10}
Statements \cite[12.2.6 and 12.2.10]{km-flips} are incorrect.
The applications of these results concern families of 3-folds with $\q$-factorial, terminal singularities. With these additional assumptions,
the proofs given there   are correct. However,  trying to formulate them with  minimal sets of  assumptions  led to  errors.

First, in  \cite[12.2.6]{km-flips} one should also assume that the fibers $X_s$ have rational singularities. This is necessary since the proof uses 
\cite[12.1.5]{km-flips}. No other changes needed.

The bigger problem is with  \cite[12.2.10]{km-flips}. During the proof  we claim to apply Kleiman's ampleness criterion to
algebraic spaces. However \cite{Kleiman66b}  states and proves
the criterion  for quasi-divisorial schemes.  This class includes all varieties that are either   projective, or proper and  $\q$-factorial, but it is not clear  that the fibers $X_s$ are quasi-divisorial.

A proof for 
$\q$-factorial, 3-dimensional, algebraic spaces is explained in \cite[5.1.3]{k-etc}, this is enough for the applications in \cite{km-flips}.
A higher dimensional Kleiman criterion needed for  \cite[12.2.10]{km-flips}
  was recently established by  Villalobos-Paz  \cite{dvp}. With this in place, the rest of the proof of \cite[12.2.10]{km-flips} is correct.
\end{say}

\begin{ack} I thank D.~Abramovich, D.~Barlet, F.~Campana, J.-P.~Demailly, O.~Fujino, S.~Kleiman, S.~Mori,  T.~Murayama, V.~Tosatti, D.~Villalobos-Paz, C.~Voisin  and C.~Xu
for  helpful comments,   corrections  and V.~Balaji for a serendipitous e-mail.
Partial  financial support    was provided  by  the NSF under grant number
DMS-1901855.
\end{ack}

 \section{Seshadri's criterion and variants}
 
   The criterion---proved by Seshadri but first published in \cite[Sec.I.7]{Ha70}---is the following.

\begin{thm}\label{sesh.crit.thm} Let $X$ be a proper algebraic space and $L$ a line bundle on $X$. Then $L$ is ample iff there is an $\epsilon>0$ such that
   $$
  \deg L|_C\geq \epsilon \cdot \mult_pC
  \eqno{(\ref{sesh.crit.thm}.1)}
   $$
  for every integral curve $C\subset X$ and every $p\in C$.
\end{thm}

   By linearity then the same holds for all 1-cycles $Z$ on $X$.
   \medskip

   An important observation is that while the criterion is frequenty stated for schemes only, it in fact holds for proper algebraic spaces, even if they are reducible and non-reduced.
     This is in marked contrast with Kleiman's criterion, which is still not known for algebraic spaces in general, though a  recent result of Villalobos-Paz  \cite{dvp}
proves  Kleiman's criterion for  $\q$-factorial algebraic spaces with log terminal singularities over a field of characteristic 0.

In algebraic geometry Seshadri's criterion is usually used to show that a line bundle $L$ is ample. Here I focus on a reformulation of Theorem~\ref{sesh.crit.thm}. 

\begin{cor}\label{sesh.crit.thm.as} Let $X$ be a proper algebraic space. Then $X$ projective  iff there is a line bundle $L$ and an $\epsilon>0$ such that
   $$
  \deg L|_C\geq \epsilon \cdot \mult_pC
  \eqno{(\ref{sesh.crit.thm.as}.1)}
   $$
  for every integral curve $C\subset X$ and every $p\in C$. \qed
\end{cor}

Now we make a slight twist and  replace the  line bundle $L$ by
 a cohomology class $\Theta\in H^2\bigl(X(\c), \q\bigr)$. The key point is that in our applications we will be able to find such a cohomology class $\Theta$, but it usually will not be a $(1,1)$ class.

 \begin{prop}\label{sesh.top.ver.prop}
   Let $X$ be a proper algebraic space over $\c$ with 1-rational singularities  (see Definition~\ref{rtl.1.rtl.defn}). Then $X$ is projective  iff there is a cohomology class $\Theta\in H^2\bigl(X(\c), \q\bigr)$ and an $\epsilon>0$ such that
   $$
   \Theta\cap [C]\geq \epsilon \cdot \mult_pC
    \eqno{(\ref{sesh.top.ver.prop}.1)}
 $$
   for every integral curve $C\subset X$ and every $p\in C$.
\end{prop}

 Here $[C]\in H_2\bigl(X(\c), \q\bigr)$ denotes the homology class of $C$, and
 we tacitly use the identification $H_0\bigl(X(\c), \q\bigr)\cong \q$ to view
 $ \Theta\cap [C]$ as a number.
The proof relies on an injection  $N_1(X,\q)\into H_2(X, \q)$ which we define next.

\begin{defn} \label{cone.of.curves}
  Let $X$ be a  proper,   complex analytic space.
  For $K=\q$ or $\r$, 
let $N_1(X,K)$ denote the $K$-vectorspace generated by
compact complex curves $C\subset X$, modulo numerical equivalence. That is, 
$$
\tsum\ a_iA_i\equiv \tsum\ b_jB_j\qtq{iff}
\tsum\ a_i\deg L|_{A_i}= \tsum\ b_j\deg L|_{B_j}
$$
for every holomorphic line bundle $L$ on $X$.

Let $H_2^{\rm alg}(X, K)\subset H_2(X, K)$ denote the vector subspace generated 
by homology classes of compact complex curves. Sending an algebraic homology class to a numerical equivalence class gives a natural surjection
$$
\res^{\rm hom}_{\rm num}: H_2^{\rm alg}(X, \q)\onto N_1(X,\q).
\eqno{(\ref{cone.of.curves}.1)}
$$
We see in Proposition~\ref{N1.to.H2.prop} that $\res^{\rm hom}_{\rm num} $ is an isomorphism if $X$ has  rational singularities. In this case the inverse map gives an  injection
$$
N_1(X,\q)\into H_2(X, \q).
\eqno{(\ref{cone.of.curves}.2)}
$$
Nore, however, that the natural map is (\ref{cone.of.curves}.1),
and (\ref{cone.of.curves}.2) exists only if $X$ has  1-rational singularities.
\end{defn}

\begin{defn} \label{rtl.1.rtl.defn} Let $X$ be a normal, complex analytic space and $g:Y\to X$ a resolution of singularities.  $X$ has {\it  rational singularities} iff $R^ig_*\o_Y=0$ for every $i>0$. The $R^ig_*\o_Y$ are independent of the choice of $Y$, so this is a property of $X$ only.
  If  $R^1g_*\o_Y=0$, then $X$ is said to have  {\it  1-rational singularities.}
  \end{defn}

 \begin{say}[Proof of Proposition~\ref{sesh.top.ver.prop}]\label{sesh.top.ver.prop.pf}
 If $C\mapsto [C]$ gives  an injection $N_1(X, \q)\into H_2\bigl(X(\c), \q\bigr)$, then we can view  $C\mapsto \Theta\cap [C]$ as a linear map
 $$
 \Theta\cap : N_1(X, \q)\to \q.
 $$
 By definition,  line bundles  span the dual space of $ N_1(X, \q)$, so
 there is a line bundle $L$ on $X$ and an $m>0$ such that
 $\deg(L|_C)=m\cdot \Theta\cap [C]$
 for every integral curve $C\subset X$. Thus
 $$
 \deg(L|_C)=m\cdot \Theta\cap [C]\geq m\epsilon \cdot \mult_pC
   $$
 for every integral curve $C\subset X$ and every $p\in C$.
 Then  $L$ is ample by Theorem~\ref{sesh.crit.thm}, so $X$ is projective.

  It remains to show that $C\mapsto [C]$ gives  an injection $N_1(X, \q)\into H_2\bigl(X(\c), \q\bigr)$ if $X$ has  1-rational singularities only. This is proved in Proposition~\ref{N1.to.H2.prop}. 
 \qed
 \end{say}

 As a consequence of Lefschetz's theorems, sending a curve $C$ to its homology class $[C]$ descends to an  injection  $N_1(X,\q)\to H_2(X, \q)$  for smooth projective varieties, see \cite[p.161]{gri-har}.
 Next we show that the same holds if $X$ has  1-rational singularities; the key ingredient is \cite[12.1.4]{km-flips}.

\begin{prop}\label{N1.to.H2.prop}
 Let $X$ be a proper, complex analytic space
that is bimeromorphic to a projective variety.
Assume that $X$ has  1-rational singularities. Then  sending a curve to its homology class gives an   injection  $N_1(X,\q)\into H_2(X, \q)$.
\end{prop}

Proof. Take a projective resolution $g:Y\to X$. Let $Z$ be a numerically trivial 1-cycle. We need to show that $[Z]\in H_2(X, \q)$ is the 0 class.

For every closed, irreducible, analytic curve
$C\subset X$ there is a closed, irreducible, analytic curve
$C_Y\subset Y$ such that $g(C_Y)=C$.  (In general $C_Y\to C$ may have degree $>1$.) Thus every $Z\in N_1(X, \q)$ lifts to $Z_Y\in N_1(Y,\q)$ such that
$g_*(Z_Y)=Z$ (as cycles).

Let $N_1(Y/X,\q)\subset N_1(X, \q)$ denote the vector subspace generated by
compact complex curves $C\subset Y$ that map to a point in $X$. 
If there is a $Z'_Y\in N_1(Y/X,\q)$ 
such that $Z_Y\equiv Z'_Y$, then 
$[Z]=g_*[Z_Y]=g_*[Z'_Y]=0$ in $H_2(X,\q)$. 

Otherwise there is a line bundle $L_Y$ on $Y$ that is
trivial on $N_1(Y/X,\q)$ but  $(L_Y\cdot Z_Y)\neq 0$. 

The key point is that, by \cite[12.1.4]{km-flips}, 
some (positive) power  $L_Y^m$ descends to  a line bundle $L$ on $X$, and then  $(L\cdot Z)=m(L_Y\cdot Z_Y)\neq 0$ gives a  contradiction. 
(This is where we use that $X$ has 1-rational singularities.) \qed

\section{Chow-Barlet spaces with marked multiplicities}

  Let $g:X\to S$ be a proper morphism. We are interested in the 
  set of proper 1-cycles $Z\subset X$ with a marked point $p\in Z$ such that $g(Z)\subset S$ is a single point and $\mult_pZ=m$.
  The corresponding coarse moduli space is usually called a
  {\it Chow variety}  in the  algebraic case
  and a {\it Barlet space} in the complex analytic setting; see 
  \cite{bar-mag} for their theory.

  We need only a rough approximation of these by a countable union of projective morphisms. This can be easily derived from the classical theory of Chow varieties for projective spaces.
  
%%   the  notation is  ${\it Chow}_1^m(X/S)$ denote the

%%   By \cite{bar-mag},  there is a complex analytic space 
%%  that is a coarse moduli space for ${\it Chow}_1^m(X/S)$; it is usually called a {\it Barlet space.}
%% However, I do not know of suitable references in general.
%% For our applications, the following much weaker version is sufficient.

\begin{prop}\label{bir.chow.prop}
  Let $g:X\to S$ be a proper morphism of complex analytic spaces  that is bimeromorphic to a projective morphism. Fix $m\in \n$. Then there are countably many diagrams of complex analytic spaces over $S$
  $$
  \begin{array}{ccc}
    C_i  & \into  &  W_i\times_SX\\
    w_i \downarrow \uparrow{\sigma_i} &&\\
    W_i&&
  \end{array}
  \eqno{(\ref{bir.chow.prop}.1)}
$$
 indexed by  $i\in I$,  such that
  \begin{enumerate}\setcounter{enumi}{1}
  \item the $w_i:C_i\to W_i$ are  proper,  of pure relative dimension 1 and flat over a dense, Zariski  open subset $W_i^\circ\subset W_i$,
  \item the fiber of $w_i$ over any $p\in W_i^\circ$ has multiplicity $m$ at $\sigma_i(p)$,
  \item the $W_i$ are irreducible, the structure maps $\pi_i:W_i\to S$ are projective, and
  \item the fibers over all the $W_i^\circ$ give all irreducible curves that have multiplicity $m$ at the marlked point.
  \end{enumerate}
  %% (Informally speaking, the union of all the $W_i^\circ$ maps surjectively onto the  relevant  Chow-Barlet space.)
\end{prop}

Proof.  By assumption there is a bimeromorphic morphism  $r:Y\to X$ such that
$Y$ is projective over $S$.
The Chow variety of curves  on $Y/S$ exists and its irreducible components are projective over $S$ (cf.\ \cite[Sec.I.5]{rc-book}). The universal curve over it parametrizes all pointed curves on $Y$.

If we have a family of pointed curves
 $$
  \begin{array}{ccc}
    C_Y  & \into  &  W\times_SY\\
    w_Y \downarrow \uparrow{\sigma_Y} &&\\
    W&&
  \end{array}
   \eqno{(\ref{bir.chow.prop}.6)}
$$
taking its image on $X$ gives
 $$
  \begin{array}{ccc}
    C  & \into  &  W\times Y\\
    w \downarrow \uparrow{s} &&\\
    W&&
  \end{array}
   \eqno{(\ref{bir.chow.prop}.7)}
  $$
   Here $w:C\to W$ is proper, of relative dimenson 1 and  flat over
   a dense, Zariski  open subset $W^\circ\subset W$. 
   The multiplicity of a fiber at a section is an upper semicontinuous function on $W^\circ$. For each $m\in \n$ let $W^m\subset W$ denote the closure of the
   set of points $p\in W^\circ$ for which $\mult_{\sigma(p)}C_p= m$.
We repeat this for all irreducible components of $W\setminus W^\circ$.
At the end we get  countably many diagrams as in  (\ref{bir.chow.prop}.1)
that satisfy (\ref{bir.chow.prop}.2--4) but not yet (\ref{bir.chow.prop}.5).

Let $X^\circ\subset X$ be the largest open set over which $r$ is an isomorphism. The above procedure gives all irreducible pointed curves that have  nonempty intersection with $X^\circ$.  Equivalently, all curves that are not contained in $X\setminus X^\circ$. (We may also get some curves contained in $X\setminus X^\circ$.)

We can now use  dimension induction to get countably many diagrams that give us all curves on  $g: (X\setminus X^\circ)\to S$. The union of these  families with the previous ones gives
 countably many diagrams
that satisfy (\ref{bir.chow.prop}.2--5). \qed

\section{Projectivity of very general fibers}

Here we prove that, under the  assumptions of Theorem~\ref{12.2.10.thm.S}, there are many projetive fibers.

\begin{prop}\label{12.2.10.vgen.prop}
  Notation and assumptions as in Theorem~\ref{12.2.10.thm.S}.
  Then there is a Euclidean    open neighborhood  $0\in U\subset S$ and
  countably many nowhere dense,  closed, analytic  subsets  $\{ H_j\subset   U: j\in J\}$,  such that $X_s$ is projective for every $s\in U\setminus \cup_j H_j$.
\end{prop}

Proof. First choose $0\in U\subset S$ such
 that $X_U$  retracts to $X_0$.
Since $X_0$ is projective,  it has an ample line bundle $L$.
Let  $\Theta\in H^2(X_U, \q)$ be the pull-back of
$c_1(L)$ to $X_U$.  Note that $\Theta$ is a topological cohomology class that is usually not the Chern class of a holomorpic line bundle.

Consider now the diagrams (\ref{bir.chow.prop}.1) indexed by the contable set $I$. Let $J\subset I$ index those diagrams for which  $H_i:=\pi_i(W_i)\subset  S$ is nowhere dense in $U$.

We aim to show that 
the restriction  $\Theta_s:=\Theta|_{X_s}$ satisfies (\ref{sesh.top.ver.prop}.1) for  $s\in U\setminus \cup_{j\in J} H_j$.

Indeed, pick a curve  $C_s\subset X_s$ and a point $p_s\in C_s$.
Set $m=\mult_{p_s}C_s$. 
By assumption  there is an $i\in I\setminus J$ and a diagram  as in (\ref{bir.chow.prop}.1)
$$
  \begin{array}{ccc}
    C_i  & \into  &  W_i\times_SX\\
    w_i \downarrow \uparrow{\sigma_i} &&\\
    W_i&&
  \end{array}
  \eqno{(\ref{12.2.10.vgen.prop}.1)}
  $$
  with  a dense, Zariski  open subset $W^\circ_i\subset W_i$, such that
  \begin{enumerate}\setcounter{enumi}{1}
  \item  $(C_s, p_s)$ is one of the fibers of $w_i$ over $W^\circ_i$,
    \item $\mult_{\sigma(p)}C_p= m$ for all $p\in W^\circ_i$, and
  \item  $\pi_i:W_i\to S$ is projective and its image contains $0\in S$.
  \end{enumerate}
 
  Thus there is a disc $D$ (say of radius $>1$) and a holomorphic map
  $\tau:D\to W_i$ such that  $\pi_i(\tau(0))=0\in S$  and
  $\pi_i(\tau (1))=s\in S$.
  After pulling back  and discarding embedded points we get
$$
  \begin{array}{ccc}
    C^D  & \into  &  D\times_SX\\
    w \downarrow \uparrow{\sigma} &&\\
    D&&
  \end{array}
  \eqno{(\ref{12.2.10.vgen.prop}.2)}
  $$
  where $w$ is flat. Let $C^D_t$  denote the fiber over $t\in D$ and
  $p_t=\sigma(t)$ the marked point. Thus $C^D_0$ lies over $0\in S$ and
  $C^D_1$ is indentified with our original curve $C_s$. 

  Note that
$$
\mult_{p_t} C^D_t=\mult_{p_1} C^D_1=\mult_{p_s}C_s \qtq{for all}t\in D^\circ.
$$
Since multiplicity is an upper semi-continous function, we conclude that
$$
\mult_{p_0}C^D_0\geq \mult_{p_t}C^D_t=\mult_{p_s}C_s.
$$
Here $C^D_0 $ is a 1-cycle on the projective scheme $X_0$,
and $\Theta_0$ is the Chern class of an ample line bundle on $X_0$.
Thus
$$
\Theta\cap [C^D_0]\geq \epsilon\cdot \mult_{p_0}C^D_0,
$$
by the easy direction of Theorem~\ref{sesh.crit.thm}, where $\epsilon$ depends only on $X_0$ and $\Theta_0$.
  Putting these together gives that
  $$
  \Theta_s\cap [C_s]=\Theta\cap [C^D_1]=\Theta\cap [C^D_0]
  \geq \epsilon\cdot \mult_{p_0}C^D_0\geq \epsilon\cdot\mult_{p_s}C_s.
  $$
  Thus $X_s$ is projective by Proposition~\ref{sesh.top.ver.prop}.
  \qed

\section{Projectivity of general fibers}

We prove a  result about the set of projective fibers of  proper analytic maps.
 To formulate it, let $g:X\to S$ be a  proper morphism of complex analytic spaces and  set
  $$
  \operatorname{PR}_S(X):=\{s\in S:  X_s \mbox{ is projective}\}.
  $$
  We prove in Proposition~\ref{RT.vb.Zd.M} that
  $\operatorname{PR}_S(X) $  either contains  a dense, Zariski open subset or
  it is contained in a  countable union of
  Zariski closed, nowhere dense subsets.
Results of this type appear in many places,  for example
  \cite{iitaka, lie-ser, ueno-83, rao-tsa-1, rao-tsa}. All the arguments below can be found in  them.

The next lemma is frequently stated in the algebraic case as in 
\cite[p.43]{gm-book}, but the proof works for proper analytic maps as well. See, for example 
\cite[p.10]{parusinski}. 

\begin{lem}\label{top.of.prop.maps.1}
Let $g:X\to S$ be a  proper morphism of complex analytic spaces.
Then there is a dense,  Zariski open subset $S^\circ\subset S$ such that
$g^\circ:X^\circ\to S^\circ$ is a topologically locally trivial fiber bundle. 
\qed
\end{lem}

Applying  this inductively we get the following.

\begin{cor}\label{top.of.prop.maps.2}
  Let $g:X\to S$ be a  proper morphism of complex analytic spaces.
  Then the sheaves $R^ig_*\z_X$ are constructible in the analytic Zariski topology. \qed
\end{cor}

  \begin{prop} \label{RT.vb.Zd.M} Let $g:X\to S$ be a  proper morphism of   normal, irreducible  analytic spaces.  
Then  there is  a dense, Zariski open subset  $S^\circ\subset S$  such that
\begin{enumerate}
\item either   $X$ is locally projective over $S^\circ$,
  \item or $\operatorname{PR}_{S}(X)\cap S^\circ$ is locally contained in a countable union of  Zariski closed, nowhere dense subsets.
\end{enumerate}
  \end{prop}

  {\it Remark \ref{RT.vb.Zd.M}.3.}   The following example clarifies  (\ref{RT.vb.Zd.M}.2).
  Let $L\subset \c^2$ be a general line and $Z$ its image in a complex torus
  $\c^2/\z^4$. Then $Z$ is irreducible and everywhere dense in the Euclidean topology. However, if $U\subset \c^2/\z^4$ is contractible, then
  $Z\cap U$ is a countable union of  Zariski closed, nowhere dense subsets of $U$.
  \medskip

  Proof. By passing to a Zariski open subset,  we may assume that 
  $R^2g_*\o_X$ is locally free and  $R^2g_*\z_X$ is  locally constant.
  By passing to the universal cover, we  may also assume that 
  $R^ig_*\z_X$ is a  constant sheaf.  The exponential sequence now shows that the holomorphic line bundles on $X$ are given by the kernel of
  $$\partial: R^2g_*\z_X\to R^2g_*\o_X,$$ and for $s\in S$,  the holomorphic line bundles on $X_s$ are given by the kernel of
  $$
  H^2(X_s, \z_{X_s})=R^2g_*\z_X|_s\to R^2g_*\o_X|_s= H^2(X_s, \o_{X_s}).
  $$
  Let $\Theta$  be a global section of $R^2g_*\z_X$. If $\partial\Theta\equiv 0$
  then there is a global line bundle $L_\Theta$ corresponding to it.
  If such an $L_\Theta$ is ample on some $X_s$, then it is ample over a
   Zariski open neighborhood of $s$ and we are in case (\ref{RT.vb.Zd.M}.1).  Assume now that this is not the case.

   If  $\partial\Theta$ is not identically $0$, then   $\partial\Theta=0$ defines a Zariski closed, nowhere dense subset $H_\Theta\subset S$.
   We claim that $\operatorname{PR}_{S}(X)$ is contained in the union of these 
  $H_\Theta$.  Indeed, pick  $s\in S$ in their complement.  Then all line bundles on $X_s$ are restrictions of a line bundle on $X$ (up to numerical equivalence) and these are not ample by assumption. \qed

\begin{complem} \label{RT.vb.Zd.M.c} Notation and assumptions as in Proposition~\ref{RT.vb.Zd.M}. Assume in addition that  $g$ is bimeromorphic to a projective morphism and we are in case (\ref{RT.vb.Zd.M}.1). Then $X$ is projective over $S^\circ$. 
\end{complem}

Proof. Choose a projective modification $\pi:X'\to X$. Set $g':=g\circ \pi$. By passing to a 
 Zariski open subset, we may assume  that 
 $R^2g_*\o_X$ and  $R^2g'_*\o_{X'}$ are  locally free, and  $R^2g_*\z_X$  and $R^2g'_*\z_{X'}$ are  locally constant.  Since $g'$ is projective, the monodromy on the algebraic classes of $R^2g'_*\z_{X'}$ is finite.
Since $X$ is normal, a
 holomorphic line bundle is trivial on $X$ iff its pull-back to $X'$ is trivial,
hence the natural map  $R^2g_*\z_X\to R^2g'_*\z_{X'}$ is injective  on  algebraic classes. Therefore 
   the monodromy on the algebraic classes of $R^2g_*\z_{X}$ is also finite.

 Once the monodromy on  the kernel of $\partial: R^2g_*\z_X\to R^2g_*\o_X$
 is finite, we can  trivialize it by a finite cover of $S$.
 We thus find a relatively ample line bundle $L_\Theta$ after a finite cover of $S$, which then gives a relatively ample line bundle on $X$. \qed

 \begin{say}[Proof of Theorem~\ref{12.2.10.thm.S}]\label{12.2.10.thm.S.pf}
   By Proposition~\ref{12.2.10.vgen.prop},  $\operatorname{PR}_{S}(X)$ contains the
   complement of a countable union of Zariski closed, nowhere dense subsets.
   Therefore, by the Baire category theorem,  $\operatorname{PR}_{S}(X)$ is not contained in a  countable union of   closed, nowhere dense subsets. Thus we are in case (\ref{RT.vb.Zd.M}.1) and $g$ is locally projective over a
   dense, Zariski open subset  $S^\circ\subset S$. Global projectivity over $S^\circ$ is given by Complement~\ref{RT.vb.Zd.M.c}.

   We finish the proof by induction applied to $S\setminus S^\circ$. \qed
   \end{say}

%\bibliography{refs}
\def\cprime{$'$} \def\cprime{$'$} \def\cprime{$'$} \def\cprime{$'$}
  \def\cprime{$'$} \def\dbar{\leavevmode\hbox to 0pt{\hskip.2ex
  \accent"16\hss}d} \def\cprime{$'$} \def\cprime{$'$}
  \def\polhk#1{\setbox0=\hbox{#1}{\ooalign{\hidewidth
  \lower1.5ex\hbox{`}\hidewidth\crcr\unhbox0}}} \def\cprime{$'$}
  \def\cprime{$'$} \def\cprime{$'$} \def\cprime{$'$}
  \def\polhk#1{\setbox0=\hbox{#1}{\ooalign{\hidewidth
  \lower1.5ex\hbox{`}\hidewidth\crcr\unhbox0}}} \def\cdprime{$''$}
  \def\cprime{$'$} \def\cprime{$'$} \def\cprime{$'$} \def\cprime{$'$}
\providecommand{\bysame}{\leavevmode\hbox to3em{\hrulefill}\thinspace}
\providecommand{\MR}{\relax\ifhmode\unskip\space\fi MR }
\providecommand{\MRhref}[2]{%
  \href{http://www.ams.org/mathscinet-getitem?mr=#1}{#2}
}
\providecommand{\href}[2]{#2}

\bigskip

  Princeton University, Princeton NJ 08544-1000, \

\email{kollar@math.princeton.edu}

\end{document}